\documentclass[13pt]{article}
\title{Braid Groups of the Sun Graph}
\author{Alice Neels and Stephen Privitera}
\date{\today}
\usepackage{amsmath,amsthm,amscd,amssymb}
\chardef\bslash=`\\ % p. 424, TeXbook
\theoremstyle{plain}

\theoremstyle{definition}
\newtheorem{defn}{Definition}[section]
\theoremstyle{remark}

\newtheorem{thm}{Theorem}[section]

\newcommand{\ra}{\rightarrow}

\newcommand{\Utop}[2]{\ensuremath{U_{#1}^{top}(#2)}}
\newcommand{\fund}[1]{\ensuremath{\pi_1(#1)}}

\swapnumbers
\theoremstyle{plain}

\begin{document}
\date{\today}
%MSC 20F36 (braid grps, artin grps)
%topics: Geometric Topology; Group Theory
%\subjclass{57 N 15,  57 R 67}
\maketitle
%%%%%%%%%%%%%%%%%%%%%%%%%%%%%%%%%%%%%%

\begin{abstract}
In this article we calculate the n-string braid groups of certain
non-contractible graphs.  We use techniques from the work of A.
Abrams, F. Connolly and M. Doig combined with Van Kampen's Theorem
to prove these results.
\end{abstract}

%%%%%%%%%%%%%%%%%%%%%%%
\section{Introduction and Statement of Results}
Recently, there has been some interest in the braid groups of a
graph $G$ (see for example \cite{abramsthesis}, \cite{abramsgeom},
\cite{ag}, \cite{condoig}, \cite{mdoig}, \cite{ghrist}).  These
groups, $B_n(G,c)$, were calculated indirectly in \cite{ghrist} in
the case where $G$ is a star. Doig critically extends this result in
\cite{mdoig} by providing \emph{explicit generators} for $B_n(G,c)$
for any star $G$.  In \cite{condoig}, Connolly and Doig provide a
method of calculating $B_n(T,c)$ where $T$ is a linear tree.  Their
method seems limited to the realm of contractible graphs (trees).

\vspace{.2cm}

In this paper we exhibit a free basis of order $n$ for $B_n(L_1,c)$
for a certain \emph{non-contractible} graph $L_1$.  We also provide
an explicit free basis of order $n+1$ for $B_2(L_n,c)$ for certain
non-contractible graphs $L_n$ with $n$ nodes.

\vspace{.2cm}

We define $L_n$ and configuration spaces and braid groups as
follows:

\begin{defn}\emph{(Sun Graph)} \\
Define $A_j=\{re^{2ij\pi/n}: 1 \leq r \leq 2\}$ for each $j$, $0\leq
j < n$. We then define the \emph{$n$ ray sun graph} as $L_n = S^1
\cup A_0 \cup A_1 \cup ... \cup A_{n-1}$.
\end{defn}

\begin{picture}(300,100)(0,0)

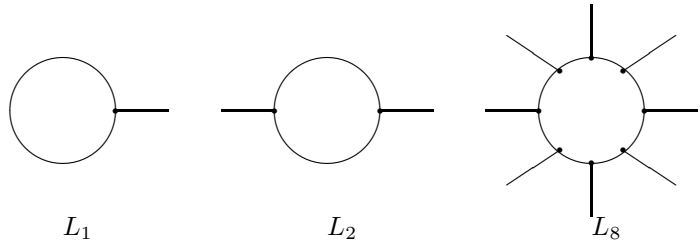
\put(90,95){\mbox{Fig. 1: Example Sun Graphs}}

%makes L_1
\put(50,50){\circle{40}}

\put(70,50){\circle*{2}}

\put(70,50){\line(1,0){20}}

\put(50,3){\mbox{$L_1$}}

%makes L_2
\put(150,50){\circle{40}}

\put(170,50){\circle*{2}}

\put(170,50){\line(1,0){20}}

\put(130,50){\circle*{2}}

\put(130,50){\line(-1,0){20}}

\put(150,3){\mbox{$L_2$}}

%makes L_8
\put(250,50){\circle{40}}

\put(270,50){\circle*{2}}

\put(270,50){\line(1,0){20}}

\put(262,65){\circle*{2}}

\put(262,65){\line(3,2){20}}

\put(250,70){\circle*{2}}

\put(250,70){\line(0,1){20}}

\put(238,65){\circle*{2}}

\put(238,65){\line(-3,2){20}}

\put(230,50){\circle*{2}}

\put(230,50){\line(-1,0){20}}

\put(262,35){\circle*{2}}

\put(262,35){\line(3,-2){20}}

\put(250,30){\circle*{2}}

\put(250,30){\line(0,-1){20}}

\put(238,35){\circle*{2}}

\put(238,35){\line(-3,-2){20}}

\put(250,3){\mbox{$L_8$}}

\end{picture}

\begin{defn} \emph{(Configuration Space and Braid Group)} \\
Let $X$ be a space. The $n$-point configuration space of $X$ is

\begin{equation*}
U_n^{top}(X) = \{c \subset X : |c|=n\},
\end{equation*}

\noindent given the quotient topology from the natural surjection

\begin{eqnarray*}
p: X^n- \Delta & \ra & U_n^{top}(X)\\
 (x_1,x_2,...x_n) & \mapsto & \{x_1,x_2,...x_n\}\\
 \end{eqnarray*}

\noindent where $\Delta =\{(x_1,x_2...x_n)\in X^n | x_i=x_j \text{
for some }i\neq j\}$ is the fat diagonal of $X^n$. The $n$-string
braid group of $X$ is

\begin{equation*}
B_n(X,c)=\pi_1(U_n^{top}(X),c)
\end{equation*}

\noindent where $c\in U_n^{top}(X)$.

\end{defn}

In order to establish our results, we will need a version of Van
Kampen's theorem to compute $\fund{A\cup B}$ when $A\cap B$ is
\emph{not} path connected.  This seems to be a well known folk
theorem, but no explicit reference is available.  For the reader's
convenience, the theorem and a proof are given in the Appendix.

This paper was funded by the NSF grant 0354132.  It is the result of
our work at the Notre Dame Research Experience for Undergraduates in
Mathematics for the summer of 2005.  We would like to thank our
advisor, Frank Connolly for his advice and support.  The authors would
also like to thank Naomi Jochnowitz for her invaluable encouragement, advice
and patience.

\section{A discussion of generators of $B_n(S_3)$}

Let $S_3$ denote the graph on 4 vertices with one vertex $v_0$ of
degree 3 and three vertices $v_1,v_2$ and $v_3$ of degree one. Label
the arms of $S_3$ as $A_j = [v_0,v_j]$ for $j=1,2,3$.

\begin{center}
\begin{picture}(100,100)(0,0)
\put(30,95){\mbox{Fig. 2: $S_3$}}
%makes a pretty S_3
\put (50,50){\line(0,-1){30}}

\put(47,55){\mbox{$v_0$}}

\put (50,50){\line(3,2){30}}

\put(47,13){\mbox{$v_1$}}

\put(52,27){\mbox{$A_1$}}

\put (50,50){\line(-3,2){30}}

\put(20,70){\mbox{$v_3$}}

\put(25,54){\mbox{$A_3$}}

\put (50,50){\circle*{2}}

\put(80,70){\mbox{$v_2$}}

\put(70,54){\mbox{$A_2$}}

\end{picture}
\end{center}

As we later make repeated use of embeddings of $S_3$ into $L_n$, it
is not necessary to specify the generators of $B_n(S_3)$.  We later
use the images of these generators to give us generators of
$B_n(L_n)$.

In \cite{condoig}, Connolly and Doig construct a one-dimensional
deformation retract, $D_n(S_3)$, of $U_n^{top}(S_3)$.  A vertex $c$
in $D_n(S_3)$ is uniquely determined by the quantities $a_i=|A_i\cap
c|$.  Therefore, we may represent each vertex by it's 3-tuple,
$(a_1,a_2,a_3)$.  We classify vertices as type I and type II
following Construction 2.2 of \cite{condoig}.  If $c$ is a type I
vertex, then $a_1+a_2+a_3=n+2$ with \mbox{$a_i<0$ $\forall$ $i$}. If
$c$ is a type II vertex, then $a_1+a_2+a_3=n$ with $a_i \geq0\forall
i$. Let $c^{(n)}$ be the configuration $(n,1,1)$.

Doig \cite{mdoig} and Connolly-Doig \cite{condoig} prove that
$B_n(S_3,c^{(n)})$ is a free group on ${n\choose 2}$ generators.  In
\cite{condoig}, Connolly and Doig construct an explicit a maximal
tree $\mathcal{T}$ of the deformation retract $D_n(S_3)$. We now
exploit the work of Connolly-Doig to give an explicit set of free
generators for $B_n(S_3,c_0)$.

The edges in $D_n(S_3)-\mathcal{T}$ are those edges beginning at the
type I vertex $(a_1,a_2,a_3)$ having $a_3>1$ and ending at the type
II vertex $(a_1-1,a_2,a_3-1)$.  We denote each such directed edge by
$e_{a_1,a_2,a_3}$, where $(a_1,a_2,a_3)$ is its associated type I
vertex. Let $\rho$ and $\rho '$ be paths in $\mathcal{T}$ from the
base vertex $(n,1,1)=c^{(n)}$ to the type I and type II vertices of
$e_{a_1,a_2,a_3}$, respectively.  Then the edge $e_{a_1,a_2,a_3}$
determines a loop, $\sigma_{a_1,a_2,a_3} = \rho\cdot
e_{a_1,a_2,a_3}\cdot \rho'^{-1}$. Because $\mathcal{T}$ is a maximal
tree in the graph $D_n(S_3)$, Van Kampen's Theorem implies that the
homotopy classes of these loops, $\sigma_{a_1,a_2,a_3}$, form a free
basis for $B_n(S_3,c^{(n)})$.

\section{The n-string braid group of the 1-ray sun graph}\label{sec4}

%fix so that D_n is in U_n by making kappa smaller

In this section we make use of a result of Abrams to give a free
basis for $B_n(L_1,c)$.  In \cite{abramsthesis}, Abrams defines the
\emph{combinatorial configuration space}, $UC_n(G)$, of a graph $G$,
and proves that this is a deformation retract of $\Utop{n}{G}$.  We
define a slightly more general space, $U_n(G)$.

\begin{defn}\emph{(Discretized Configuration Space)}\\
Let $G$ be a graph with all edges of length $\kappa$.  Then the
\emph{discretized configuration space}, $U_n(G)$, is defined as the
subset of $\Utop{n}{G}$ consisting of all configurations in which
for any two points in the configuration we can find an open edge $e$
between them.
\end{defn}

When the distance between any two essential vertices on $G$ is at
least $\kappa (n+1)$ and the length of any cycle in $G$ is at least
$\kappa (n+1)$, Abrams (\cite{abramsthesis}) proves that $UC_n(G)$
is a deformation retract of $\Utop{n}{G}$. As we change only the
standard edge length, this proof also shows that $U_n(G)$ is a
deformation retract of $\Utop{n}{G}$.

Set $\kappa=\pi/(2n-1)$. Let $L= S^1 \cup [1,2\kappa(n-1)+1]$, a
subset of $\mathbb{C}$. Clearly, \mbox{$L\approx L_1$} and therefore
$B_n(L,c)\cong B_n(L_1,c')$ for any $c\in \Utop{n}{L}$ and $c'\in
\Utop{n}{L_1}$. We make $L$ into a graph by dividing $S^1$ into
$4n-2$ edges of length $\kappa$ such that 1 is a vertex.  We divide
the ray into $2n-2$ edges of length $\kappa$.

\begin{center}
\begin{picture}(200,100)(0,0)
\put(50,95){\mbox{Fig. 3: Diagram of $L$}}
%makes a pretty S_3
\put (100,50){\circle{40}}

\put (120,50){\line(1,0){30}}

\put (120,50){\circle*{2}}

\put (130,50){\circle*{2}}

\put (140,50){\circle*{2}}

\put (150,50){\circle*{2}}

\put (80,50){\circle*{2}}

\put(112,66){\circle*{2}}

\put(100,70){\circle*{2}}

\put(88,66){\circle*{2}}

\put(80,50){\circle*{2}}

\put(112,34){\circle*{2}}

\put(100,30){\circle*{2}}

\put(88,34){\circle*{2}}

\put (60,45){\mbox{\LARGE{$I \{$}}}

\end{picture}
\end{center}

We define $I\subset S^1$ to be the closed interval of length
$\kappa/2$ centered at -1.  We define $Y=Cl_L(L-I)$.  Let $S_3$ be
as before. However, specify that each arm of $S_3$ has length
$2\kappa(n-1)$. There is a natural injective, distance preserving
map $i:S_3\ra Y$ such that $i(v_0)=1$ and
$i(A_1)=[1,2\kappa(n-1)+1]$. The map $i$ induces an injective map
$i_c: \Utop{n}{S_3} \ra \Utop{n}{L}$. We define $c^n=i_c(c^{(n)})$.
Note that $i_c(D_n(S_3))\subset U_n(L)$, where the unit length for
$D_n(S_3)$ is taken to be $2\kappa$ (instead of 1 as in
\cite{mdoig}).

Let $\beta_0$ be the homotopy class of the loop in $U_n(L)$ at the
base point $c^n$ in which the point at 1 moves once counterclockwise
around $S^1$ and the other points remain fixed. The map $i_c$
induces a homomorphism \mbox{$i_*: B_n(S_3,c^{(n)})\ra B_n(L,c^n)$}
of fundamental groups. Let $\beta_i = i_*([\sigma_{i,n-i,2}])$ for
$1\leq i < n$, where $[\sigma_{i,n-i,2}]$ is a generator of
$B_n(S_3,c^{(n)})$, as described in section 2.

The purpose of this section is to prove the following:

\begin{thm}
 $B_n(L, c^n) = F(\beta_0,\beta_1,...\beta_{n-1})$,
the free group on the $n$ letters $\beta_0,\beta_1,...\beta_{n-1}$.
\end{thm}

\begin{proof} Define $\ast_0=\{x\}\cup c^{n-1}$ and $\ast_1=\{y\}\cup
c^{n-1}$ where $x$ and $y$ are the uppermost and lowermost points in
I. Define $\mathcal{A}$ and $\mathcal{B}$ as follows:

\begin{eqnarray*}
\mathcal{A} & = & \{c\in U_n(L)| c\cap I \neq \emptyset\}\\
\mathcal{B} & = & U_n(L)\cap \Utop{n}{Y} \text{ where } Y=Cl_L(L-I)\\
\end{eqnarray*}

Notice that $\mathcal{B}$ contains all configurations in $U_n(L)$
containing $x$ or $y$, and that at most one point can be in $I$ in
any given configuration in $U_n(L)$ because the length of $I$ is
less than $\kappa$. Note that $U_n(L)=\mathcal{A}\cup\mathcal{B}$.
We will apply the Generalized Van Kampen's Theorem (see Appendix) to
these sets.

We have $\mathcal{A}\cap\mathcal{B}=\mathcal{C}\sqcup \mathcal{D}$,
where:

\begin{eqnarray*}
\mathcal{C} & = & \{c\in U_n(L)| x\in c\}\\
\mathcal{D} & = & \{c\in U_n(L)| y\in c\}\\
\end{eqnarray*}

Note that both $\mathcal{C}$ and $\mathcal{D}$ are deformation
retracts of $\mathcal{A}$.

We note that $i_c(D_n(S_3))$ is a deformation retract of
$\mathcal{B}$. Therefore \mbox{$\pi_1(\mathcal{B},\ast_0)\cong F_{n
\choose 2}$} by the formula given in \cite{mdoig} for the $n$-string
braid group of a 3-star.

Let $\tau$ be a path from $\ast_0$ to $c^{(n)}$.  Define
$\beta_{i,j,k}$ to be the homotopy class
$[\tau\sigma_{i,j,k}\tau^{-1}]$.  The generators for
$\fund{\mathcal{B},\ast_0}$ are $\beta_{i,j,k}$ where $i+j+k=n+2$
and $k>1$.

There a natural map $k: D_{n-1}(S_3)\ra\mathcal{A}$ given by
$c\mapsto i_c(c)\cup\{x\}$.  This map is a homotopy equivalence.
Therefore \mbox{$\pi_1(\mathcal{A},\ast_0)\cong
B_{n-1}(S_3,c^{(n-1)})$}. By the formula given in Doig for the
$n-1$-string braid group of a 3-star,
\mbox{$\pi_1(\mathcal{A},\ast_0)\cong F_{n-1 \choose 2}$} . The
generators for $\pi_1(\mathcal{A},\ast_0)$ are
$\alpha_{i,j,k}=k_\ast([\sigma_{i,j,k}])$ where $i+j+k=n+1$ and
$k>1$.

Since $\mathcal{C}$ and $\mathcal{D}$ are deformation retracts of
$\mathcal{A}$, we also have $\fund{\mathcal{C},\ast_0}\cong F_{n-1
\choose 2}$ and $\fund{\mathcal{D},\ast_1}\cong F_{n-1 \choose 2}$.
 The generators $\gamma_{i,j,k}$ for $\pi_1(\mathcal{C},\ast_0)$ and
the generators $\delta_{i,j,k}$ for $\pi_1(\mathcal{D},\ast_1)$ are
defined similarly to those of $\mathcal{A}$.

Next, we must find paths $t_\mathcal{A}$ from $\ast_0$ to $\ast_1$
through $\mathcal{A}$ and $t_\mathcal{B}$ from $\ast_0$ to $\ast_1$
through $\mathcal{B}$.  These paths lead us to the maps

\begin{eqnarray*}
j_0^\mathcal{A}:\fund{\mathcal{C},\ast_0} & \ra & \fund{\mathcal{A},\ast_0}\\
j_1^\mathcal{A}:\fund{\mathcal{D},\ast_1} & \ra & \fund{\mathcal{A},\ast_0}\\
j_0^\mathcal{B}:\fund{\mathcal{C},\ast_0} & \ra & \fund{\mathcal{B},\ast_0}\\
j_1^\mathcal{B}:\fund{\mathcal{D},\ast_1} & \ra & \fund{\mathcal{B},\ast_0}\\
\end{eqnarray*}

Explicitly,

\begin{eqnarray*}
j_0^\mathcal{A}([\sigma]) & = & [\sigma]\\
j_1^\mathcal{A}([\tau]) & = & [t_\mathcal{A}^{-1}\tau t_\mathcal{A}]\\
j_0^\mathcal{B}([\sigma]) & = & [\sigma]\\
j_1^\mathcal{B}([\tau]) & = & [t_\mathcal{B}^{-1}\tau t_\mathcal{B}]\\
\end{eqnarray*}

Notice that the first two maps are isomorphisms.  Clearly, because
$\mathcal{C}$ and $\mathcal{D}$ are deformation retracts of
$\mathcal{A}$,

\begin{eqnarray*}
j_0^\mathcal{A}(\gamma_{i,j,k}) & = & \alpha_{i,j,k}\\
j_1^\mathcal{A}(\delta_{i,j,k}) & = & \alpha_{i,j,k}\\
\end{eqnarray*}

It follows from the definition of $\beta_{i,j,k}$ that

\begin{eqnarray*}
j_0^\mathcal{B}(\gamma_{i,j,k}) & = & \beta_{i,j,k+1}\\
j_1^\mathcal{B}(\delta_{i,j,k}) & = & \beta_{i,j+1,k}\\
\end{eqnarray*}

We apply Generalized Van Kampen's Theorem now to
$U_n(L)=\mathcal{A}\cup\mathcal{B}$. We get

\begin{equation*}
B_n(L,\ast_0)\cong\frac{\fund{\mathcal{A},\ast_0}\ast\fund{\mathcal{B},\ast_0}\ast
F(t)}{N},
\end{equation*}

\noindent where $N$ is the smallest normal group containing the
${n-1 \choose 2}$ relations of the form

\begin{equation*}
\alpha_{i,j,k}\beta_{i,j,k+1}^{-1}
\end{equation*}

\noindent and the ${n-1 \choose 2}$ relations of the form

\begin{equation*}
\alpha_{i,j,k}t\beta_{i,j+1,k}^{-1}t^{-1}.
\end{equation*}

We will first eliminate all $\alpha_{i,j,k}$. For a given
$\alpha_{i,j,k}$, we have two relations involving $\alpha_{i,j,k}$.
They are:

\begin{equation*}
\alpha_{i,j,k}\beta_{i,j,k+1}^{-1}
\end{equation*}

\noindent and

\begin{equation*}
\alpha_{i,j,k}t\beta_{i,j+1,k}^{-1}t^{-1}
\end{equation*}

So we may eliminate each $\alpha_{i,j,k}$.  We are left with ${n-1
\choose 2}$ relations of the form:

\begin{equation*}
\beta_{i,j,k+1}t\beta_{i,j+1,k}^{-1}t^{-1}
\end{equation*}

Now eliminate the generators of the form $\beta_{i,j,k+1}$. This
will eliminate exactly ${n-1 \choose 2}$ generators because
$k\geq2$.

We are left with $n-1$ generators of the form $\beta_{i,j,2}$ where
$i+j=n$ and $i,j\geq 1$.  Therefore $B_n(L, \ast_0)$ has a free
basis consisting of the $n$ elements in $\{\beta_{i,j,2}|i+j=n$ and
$i,j\geq 1\}\cup\{t\}$.

By a simple change of base point from $\ast_0$ to $c^n$, we find
that $B_n(L, c^n)=F(t',\beta'_{1,n-1,2},...\beta'_{n-1,1,2})$ where
$t' = [\tau t\tau^{-1}]$ and $\beta'_{i,n-i,2} =[\tau
\beta_{i,n-i,2} \tau^{-1}]$ where $\tau$ is some path from $c^{n}$
to $\ast_0$.

Clearly, these are exactly the generators specified in our statement
of the theorem.\end{proof}

\section{The 2-Point Braid Group of the n-ray sun graph}\label{sec5}

Let $L_n$ be as before.  Let $\zeta = e^{2\pi i /n}$.
For each $j$, $0\leq j < n$, define
\begin{center}  $c_j = \{ \zeta^j, 2\zeta^j \}$\\
        $L^j = S^1 \cup A_j$\\
        $L^{-j} = Cl_{L_n}(L_n - A_j)$. \end{center}

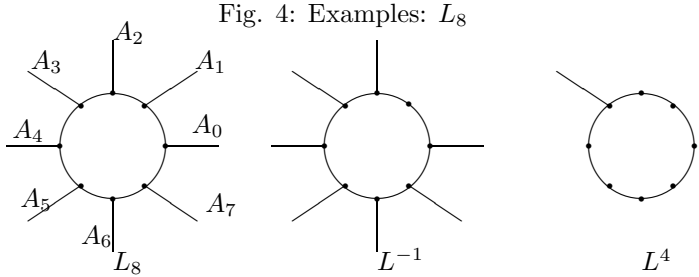
\vspace{.5cm}

\begin{picture}(300,100)(0,0)

\put(90,97){\mbox{Fig. 4: Examples: $L_8$}}

%makes L_8
\put(50,50){\circle{40}}

\put(70,50){\circle*{2}}

\put(70,50){\line(1,0){20}}

\put(80,55){\mbox{$A_0$}}

\put(62,65){\circle*{2}}

\put(62,65){\line(3,2){20}}

\put(81,79){\mbox{$A_1$}}

\put(50,70){\circle*{2}}

\put(50,70){\line(0,1){20}}

\put(50,90){\mbox{$A_2$}}

\put(38,65){\circle*{2}}

\put(38,65){\line(-3,2){20}}

\put(19,79){\mbox{$A_3$}}

\put(30,50){\circle*{2}}

\put(30,50){\line(-1,0){20}}

\put(12,53){\mbox{$A_4$}}

\put(62,35){\circle*{2}}

\put(62,35){\line(3,-2){20}}

\put(15,27){\mbox{$A_5$}}

\put(50,30){\circle*{2}}

\put(50,30){\line(0,-1){20}}

\put(38,13){\mbox{$A_6$}}

\put(38,35){\circle*{2}}

\put(38,35){\line(-3,-2){20}}

\put(85,25){\mbox{$A_7$}}

\put(50,3){\mbox{$L_8$}}

%makes L-1
\put(150,50){\circle{40}}

\put(170,50){\circle*{2}}

\put(170,50){\line(1,0){20}}

\put(162,66){\circle*{2}}

\put(150,70){\circle*{2}}

\put(150,70){\line(0,1){20}}

\put(138,65){\circle*{2}}

\put(138,65){\line(-3,2){20}}

\put(130,50){\circle*{2}}

\put(130,50){\line(-1,0){20}}

\put(162,35){\circle*{2}}

\put(162,35){\line(3,-2){20}}

\put(150,30){\circle*{2}}

\put(150,30){\line(0,-1){20}}

\put(138,35){\circle*{2}}

\put(138,35){\line(-3,-2){20}}

\put(150,3){\mbox{$L^{-1}$}}

%makes L^4
\put(250,50){\circle{40}}

\put(270,50){\circle*{2}}

\put(262,65){\circle*{2}}

\put(250,70){\circle*{2}}

\put(238,65){\circle*{2}}

\put(238,65){\line(-3,2){20}}

\put(230,50){\circle*{2}}

\put(262,35){\circle*{2}}

\put(250,30){\circle*{2}}

\put(238,35){\circle*{2}}

\put(250,3){\mbox{$L^4$}}

\end{picture}

\vspace{.5cm}

\noindent Let $t$ denote the loop given by $t(s) = \{e^{2\pi i
s},2\}$ for $0\leq s \leq 2\pi$. Finally, define $c_*=\{1,\zeta\}$.

For each $j$, $0\leq j < n$, there is an injection $i_j:S_3 \ra L_n$
such that
\begin{center}  $i_j(v_0)=\zeta^j$\\
        $i_j([v_0,v_1])=A_j$\\
        $i_j([v_0,v_2])=[\zeta^j, \zeta^{j+1}]$\\
        $i_j([v_0,v_3])= [\zeta^j, \zeta^{j-1}]$ \end{center}
This map induces a map of configuration spaces
$i_j':\Utop{n}{S_3} \ra  \Utop{n}{L_n}$.
From the above discussion, $B_2(S_3,c^{(2)})=F(\beta_{1,1,2})$.
Let $\alpha_j = i_j'(\sigma_{1,1,2})$, where $\sigma_{1,1,2}$ is
the loop in $\Utop{n}{S_3}$ described in section 2.
We have the following theorem. \\

\begin{thm}
 Let $\tau_j$ be the obvious (counterclockwise) path
from $c_0$ to $c_j$. Then the elements of the form $[\tau_j \alpha_j
\tau_j^{-1}]$ together with $[t]$ form a free basis for the group
$B_2(L_n,c_0)$.
\end{thm}

\begin{proof} The base case when $n=1$ has already been established in
the previous section. Now suppose that the theorem is true for all
$L_k$ with $k \leq n$. For the time being, we will work with the
base configuration $c_*$. Define the following two sets:

\begin{center}  $A = U_2^{top}(L^{-0}) \cup U_2^{top}(L^{-1})$ \\
        $B = \{c\in U_2^{top}(L_n)| c\cap  A_0 \neq \emptyset$ and
            $c\cap  A_1 \neq \emptyset \}$
\end{center}

The sets $B$ and $A\cap B$ are contractible subsets of
$U_2^{top}(L_n)$. In fact, $B$ is homeomorphic to $[0,1]\times [0,1]$,
and $A\cap B$ is homeomorphic to $[0,1]$.
By Van Kampen's theorem, it follows that $B_2(L_n,c_*)\cong \pi_1(A,c_*)$.
Since $A$ is path-connected, it is enough to compute $\pi_1(A,c_0)$.\\
Define the following sets:

\begin{eqnarray*}
D & = & U_2^{top}(L^{-0})\\
E & = & U_2^{top}(L^{-1}).\\
\end{eqnarray*}

By the induction hypothesis applied to $L_{n-1}$, we know the free
generators for $B_2(L_{n-1},c_0)$ are of the form  $[\tau_j \alpha_j
\tau_j^{-1}]$ for $0\leq j \leq n -1$ together with $[t]$.  There is
a homeomorphism of $L_{n-1}$ onto $L^{-0}$ which maps $A_j\subset
L_{n-1}$ to $A_{j+1}\subset L^{-0}$ for all $j$. This homeomorphism
induces a homeomorphism of $U_2^{top}(L_{n-1})$ onto $D$, which in
turn induces an isomorphism of fundamental groups
$B_2(L_{n-1},c_0)\ra B_2(L^{-0},c_0)$. The isomorphism maps the
generator $[\tau_j \alpha_j \tau_j^{-1}]$ in $B_2(L_{n-1},c_0)$ to
$\delta_{j+1} = [\tau_{j+1} \alpha_{j+1} \tau_{j+1}^{-1}]$ in
$B_2(L^{-0},c_0)$.  Similarly, there is a homeomorphism of $L_{n-1}$
onto $L^{-1}$ which maps $A_j\subset L_{n-1}$ to $A_{j+1}\subset
L^{-1}$ for $1 \leq j < n$ and maps $A_0 \subset L_{n-1}$ to
$A_{0}\subset L^{-1}$. This homeomorphism induces a homeomorphism of
$U_2^{top}(L_{n-1})$ onto $E$, which in turn induces an isomorphism
of fundamental groups $B_2(L_{n-1},c_0)\ra B_2(L^{-1},c_0)$. The
isomorphism maps the generator $[\tau_j \alpha_j \tau_j^{-1}]$ in
$B_2(L_{n-1},c_0)$ to $\epsilon_{j+1} =[\tau_{j+1} \alpha_{j+1}
\tau_{j+1}^{-1}]$ for $j>0$ and to $\epsilon_0 = [\tau_0 \alpha_0
\tau_0^{-1}]$ for $j=0$.  Let the image of $[t]$ under these
isomorphisms be denoted by $\delta\in B_2(L^{-0},c_0)$ and
$\epsilon \in B_2(L^{-1},c_0)$.\\

Now consider $D\cap E$.  There is a homeomorphism between this space
and $U_2^{top}(L_{n-2})$ which maps $A_j\subset L_{n-1}$ to
$A_{j+2}\subset L^{-0}$ for all $j$. The induced isomorphism of
fundamental groups will send $[\tau_j \alpha_j \tau_j^{-1}]$ in
$B_2(L_{n-2},c_0)$ to $\gamma_{j+2} = [\tau_{j+2} \alpha_{j+2}
\tau_{j+2}^{-1}]$ in $B_2(L^0\cap L^1,c_0)$. Let $\gamma$ denote the
image of $[t]$ under this isomorphism.
Since\\

\begin{eqnarray*}
i^D(\gamma_j)& = & \delta_j\\
i^E(\gamma_j) & = & \epsilon_j\\
i^D(\gamma) & =& \delta\\
i^E(\gamma) & = & \epsilon,\\
\end{eqnarray*}

it follows by (Classical) Van Kampen's Theorem that

$B_2(L_n,c_0)= F(\gamma_0, \delta_1, \delta_2, ..., \delta_{n-1}, t)$.
\end{proof}

\section{Appendix}

\begin{thm} \emph{Generalized Van Kampen's Theorem:}

Let $X$ be a polyhedron.  Let $A$ and $B$ be subpolyhedra of $X$
such that $A\cup B=X$ and $A\cap B=\sqcup^n_{i=0}C_i$ where each
$C_i$ is path connected and nonempty.  Choose a base point $c_i\in
C_i$ for each $C_i$.  Find paths $t_1^A,...,t_n^A,t_1^B,...,t_n^B$
such that $t_i^A$ goes from $c_0$ to $c_i$ in $A$ and $t_i^B$ goes
from $c_0$ to $c_i$ in $B$. Let $t_0^A=t_0^B$ be the constant path.
The following maps of fundamental groups are induced by the
inclusion of $C_i$ into $A$ and $B$ followed by a base point change:

\begin{equation*}
\begin{array}{rrcl}
j_i^A: & \fund{C_i,c_i} & \ra & \fund{A,c_0}\\
& [\sigma] & \mapsto & [(t_i^A)^{-1}\sigma t_i^A]\\
j_i^B: & \fund{C_i,c_i} & \ra & \fund{B,c_0}\\
& [\sigma] & \mapsto & [(t_i^B)^{-1}\sigma t_i^B]\\
\end{array}
\end{equation*}

Then there exists an isomorphism

\begin{equation*}
\Phi:\frac{\fund{A,c_0}\ast \fund{B,c_0}\ast
F(t_1,...t_n)}{N}\ra\fund{X,c_0}
\end{equation*}

\noindent where the $t_i$'s are indeterminates and $N$ is the
smallest normal subgroup containing all words of the form
$j_i^A([\sigma])(t_i j_i^B([\sigma])t_i^{-1})^{-1}$ for $0\leq i\leq
n$.  $t_0$ is defined as the identity element of $\fund{A,c_0}\ast
\fund{B,c_0}$.

The map $\Phi$ is specified by

\begin{eqnarray*}
\Phi(\alpha) & = & i^A\alpha \text{ for all } \alpha\in\fund{A,c_0}\\
\Phi(\beta) & = & i^B\beta \text{ for all } \beta\in\fund{B,c_0}\\
\Phi(t_i)& = & [t_i^A(t_i^B)^{-1}] \text{ for } 1\leq i\leq n\\
\end{eqnarray*}

\noindent where $i^A$ and $i^B$ are the homomorphisms induced by the
inclusions of $A$ and $B$ into $X$.
\end{thm}

\begin{proof}

Let $X$ be a polyhedron with subpolyhedra $A$ and $B$ such that
$X=A\cup B$ and $A\cap B=\sqcup^n_{i=0}C_i$ where each $C_i$ is path
connected and nonempty.  Let $c_i\in C_i$ be base points.

Recall that an arc is a homeomorphic copy of the unit interval.  Let
$t_i^A\subset A$ be arcs connecting $c_0$ to $c_i$ in $A$ and
$t_i^B\subset B$ be arcs connecting $c_0$ to $c_i$ in $B$ for $1\leq
i\leq n$ such that $t_i^A\cup t_i^B\approx S^1$.  We also require
that $A\cap (t_i^A\cup t_i^B)$ and $B\cap (t_i^A\cup t_i^B)$ are
arcs.

We will work by induction.  Assume that the theorem holds for some
$n\geq 1$.  The base case ($n=0$) is obtained from Van Kampen's
Theorem.

Set

\begin{eqnarray*}
A' & = & A\cup t_1^B\\
B' & = & B\cup t_1^A\\
\end{eqnarray*}

Now $A'\cup B'=X$, but $A'\cap B'= C_1'\cup C_2\cup...\cup C_n$
where $C_1'=C_0\cup C_1 \cup t_1^A \cup t_1^B$.  Note that $A'\cap
B'$ has one less component that $A\cap B$.

We know that

\begin{eqnarray*}
\fund{A',c_0} & = & \fund{A,c_0}\ast F(t)\\
\fund{B',c_0} & = & \fund{B,c_0}\ast F(t)\\
\fund{C_1',c_0} & = & \fund{C_0,c_0}\ast \fund{C_1\cup t_1^A\cup t_1^B,c_0}\\
\end{eqnarray*}

Then we apply the inductive hypothesis to $X$ as $A'\cup B'$.  The
theorem follows.
\end{proof}

\vspace{2cm}
 Alice Neels is and undergraduate at Reed College and may be reached
 at neelsa@reed.edu.  Stephen Privitera is an undergraduate at the
 University of Rochester and may be reached at
 sprivite@mail.rochester.edu.

\end{document}